\documentclass{amsart}
\usepackage{amsfonts}
\usepackage{amsmath}
\usepackage{enumerate}
\usepackage{amsthm}
\usepackage{hyperref}
\usepackage{cite}

\begin{document}

\newtheorem{definition}{Definition}[section]
\newtheorem{theorem}[definition]{Theorem}
\newtheorem{proposition}[definition]{Proposition}
\newtheorem{remark}[definition]{Remark}
\newtheorem{lemma}[definition]{Lemma}
\newtheorem{corollary}[definition]{Corollary}
\newtheorem{example}[definition]{Example}

\numberwithin{equation}{section}

\title[A note on Kazdan-Warner equations...]{A note on Kazdan-Warner type equations on compact Riemannian manifolds }
\author[W. Yu]{Weike Yu}

%\thanks{The author is supported by Jiangsu Funding Program for Excellent Postdoctoral Talent (2022ZB281), and the Fundamental Research Funds for the Central Universities (30922010410).}
\date{}

\begin{abstract}
In this note, we prove an existence result for generalized Kazdan-Warner equations on compact Riemannian manifolds by using the flow approach or the upper and lower solution method. In addition, we give a priori estimates for this type equations.

\end{abstract}
\maketitle
\section{Introduction}
In \cite{[KW]}, Kazdan-Warner initiated the study of the following semi-linear elliptic partial differential equation which is now called the Kazdan-Warner equation
\begin{align}\label{1.1}
\Delta u-c-h e^{u}=0
\end{align}
on a compact Riemannian surface $(M^2, g)$, where $\Delta$ is the Laplace-Beltrami operator with respect to $g$, and $c, h$ are two given smooth functions on $M$. Their motivation of studying this equation was to solve the prescribed Gaussian curvature problem on Riemannian surfaces, and they had given satisfying characterizations to the solvability of this equation on Riemannian surfaces. Later, the Kazdan-Warner type equation was also investigated in various contexts, such as \cite{[Ber], [Bra], [DJLW1], [Doa], [JT], [LY], [Mi], [Mor], [MP], [Tau], [Wan], [WZ]}, and many others.

Recently, Angella-Calamai-Spotti \cite{[ACS]} studied a Hermitian analogue of the classical Yamabe problem on compact Hermitian manifolds which is called the Chern-Yamabe problem, and later the prescribed Chern scalar curvature problem was also investigated by several authors (cf. \cite{[Fus], [CZ], [Ho1], [Ho2], [LZZ], [LM], [Yu]}, etc.).  According to \cite{[ACS]}, we know that solving these problems is equivalent to proving the existence of the following type equation on a compact Hermitian manifold $(M^n, \omega)$:
\begin{align}\label{1.2}
\Delta^{Ch} u-c-h e^{u}=0,
\end{align}
where $\Delta^{Ch}$ denotes the Chern Laplacian with respect to the Hermitian metric $\omega$. We remark that $\Delta^{Ch}u=\Delta_du-\langle du, \theta_\omega\rangle$ for any $u\in C^2(M)$, where $\Delta_d$ is the Hodge-de Rham Laplacian, $\theta_\omega$ is the Lee form (or torsion $1$-form) associated with $\omega$, and $\langle\cdot, \cdot\rangle$ denotes the inner product on $1$-form induced by $\omega$.

In this paper, we consider the following generalized Kazdan-Warner equation on a compact Riemannian manifold $(M^n, g)$ (which was introduced by Bryan and Wentworth in \cite{[BW]} to construct the moduli spaces for K\"ahler surfaces).
\begin{align}\label{1.3}
 \Delta u-\langle du, \theta\rangle-S-A e^{\alpha u}+B e^{-\beta u}=0,
\end{align}
where $\theta$ is a $1$-form on $M$ with $d^*\theta=0$, $\alpha, \beta$ are two positive numbers, and $A, B, S$ are three smooth functions on $M$. Note that a well-known example of the operator $L=\Delta -\langle d\cdot, \theta\rangle$ with $d^*\theta=0$ is the Chern Laplacian $\Delta^{Ch}$ with respect to a Gauduchon metric $\omega$ (since its corresponding Lee form $\theta_\omega$ satisfies $d^*\theta_\omega=0$). Moreover, when $B\not\equiv 0$, the above equation \eqref{1.3} has somewhat different characteristics from the classical Kazdan-Warner equation \eqref{1.1} and Chern-Yamabe equation \eqref{1.2}. 

In terms of the flow approach or the upper and lower solution method, we prove the following existence result 
\begin{theorem}\label{theorem1.1}
Let $(M,g)$ be a compact Riemannian manifold, and let $A, B, S$ be smooth functions on $M$ satisfying
\begin{align}\label{1.4}
A\geq 0\ (\not\equiv 0),\ B\geq 0,\ \int_M S<0.
\end{align}
Then the generalized Kazdan-Warner equation \eqref{1.3} has a unique smooth solution on $M$.
\end{theorem}
\begin{remark}
\begin{enumerate}
\item The proof is partly inspired by our recent work on prescribed Chern scalar curvature problem \cite{[Yu]}.
\item In \cite{[BW]}, Bryan-Wentworth proved that \eqref{1.3} with $\theta\equiv 0$ has a unique smooth solution on $M$ provided that $\eqref{1.4}$ and $\int_M(A-B)>0$ hold.
\end{enumerate}
\end{remark}
In addition, we prove the following priori estimates for the equation \eqref{1.3}.

\begin{proposition}\label{prop1.2}
Let $(M, g)$ be a compact Riemannian manifold, and let $A, B, S$ be smooth functions on $M$ such that
\begin{align}
A>0,\ B\geq 0,\ \int_M S<0.
\end{align}
Then there exists a constant $C_k$ such that for any function $u\in C^\infty(M)$ satisfying \eqref{1.3}, we have
\begin{align}\label{3.2}
\|u\|_{C^k(M)}\leq C_k
\end{align}
for any $k=0,1,2,\cdots$.
\end{proposition}

\textbf{Acknowledgements.} The author would like to thank Prof. Yuxin Dong and Prof. Xi Zhang for their continued support and encouragement.

\section{Proof of Theorem \ref{theorem1.1}} 
In this section, we will give the proof of Theorem \ref{theorem1.1}. Firstly, we consider the following flow
\begin{equation}\label{2.1}
\left\{
\begin{aligned}
&\frac{\partial u}{\partial t}= \Delta u-\langle du, \theta\rangle-S-A e^{\alpha u}+B e^{-\beta u},\\
&u(x,0)=u_0,
\end{aligned}
\right.
\end{equation}
where $d^*\theta=0$, $\alpha, \beta$ are two positive numbers, $A, B, S\in C^\infty(M)$ satisfy \eqref{1.4}, and $u_0$ is an arbitrary smooth function on $M$. According to the standard parabolic equation theory, the flow \eqref{2.1} has a smooth solution $u: M\times [0,T)\rightarrow \mathbb{R}$ for some $T>0$. Here we recall the maximum principle for parabolic equations, which will be used later.
\begin{lemma}[{\cite[Theorem 3.1.1, page 44]{[Top]}}]\label{lemma2.1}
Suppose $\{g(t)\}_{t\in [0,t_0]}\ (0<t_0<\infty)$ be a smooth family of Riemannian metrics, and $X(t)$ is a smooth family of vector fields on a compact manifold $M$. Assume that $u\in C^\infty(M\times[0,t_0])$ satisfies
\begin{align}
\frac{\partial u}{\partial t}\leq \Delta_{g(t)} u+g(t)(X(t), \nabla u)+F(u,t),
\end{align}
where $F: \mathbb{R}\times [0,t_0]\rightarrow \mathbb{R}$ is a smooth function. Then if $u(\cdot,0)\leq c_0\in \mathbb{R}$, we have $u(\cdot,t)\leq \phi(t)$ for any $t\in [0,t_0]$, where $\phi(t)$ is the solution of 
\begin{equation}
\left\{
\begin{aligned}
&\frac{d\phi}{dt}=F(\phi(t),t),\\
&\phi(0)=c_0.
\end{aligned}
\right.
\end{equation}
\end{lemma}
Making use of the above maximum principle, we prove that the solution of \eqref{2.1} is global (namely, $T=+\infty$), and is of a uniform lower bound with respect to $t\in [0,+\infty)$.
\begin{lemma}
The flow \eqref{2.1} has a unique smooth solution $u(x,t)$ defined on $M\times [0,+\infty)$.
\end{lemma}
\proof Set $L=\Delta -\langle d\cdot, \theta\rangle$. Since $A,B\geq 0, \alpha,\beta>0$, we deduce from \eqref{2.1} that
\begin{align}\label{2.3}
\left(\frac{\partial}{\partial t}-L\right)\left|\frac{\partial u}{\partial t}\right|^2(x,t)=-2\left(\alpha Ae^{\alpha u}+\beta B e^{-\beta u}\right)\left|\frac{\partial u}{\partial t}\right|^2-\left|\nabla \frac{\partial u}{\partial t}\right|^2\leq0
\end{align}
for any $(x,t)\in M\times [0,T)$, and
\begin{align}\label{3.3..}
\left|\frac{\partial u}{\partial t}\right|^2(x,0)=\left|L u_0-S-A e^{\alpha u_0}+B e^{-\beta u_0}\right|^2\leq C_1^2,
\end{align}
where $C_1=\sup_M(\left|L u_0-S-A e^{\alpha u_0}+B e^{-\beta u_0}\right|)<\infty$. Applying Lemma \ref{lemma2.1} to \eqref{2.3} and \eqref{3.3..}, we have 
\begin{align}\label{2.5}
\sup_{M\times[0,T)}\left|\frac{\partial u}{\partial t}\right|\leq C_1.
\end{align}
Hence,
\begin{align}\label{2.6}
|u(x,t)|\leq |u(x,t)-u(x,0)|+|u(x,0)|\leq \int^t_0\left |\frac{\partial u}{\partial t}\right|dt+\sup_M|u_0|\leq C_2,
\end{align}
where $C_2=C_1T+\sup_M|u_0|$. By \eqref{2.1}, \eqref{2.5} and \eqref{2.6}, we get
\begin{align}\label{2.8}
|L u|(x,t)\leq \sup_M|S|+(\sup_M A) e^{\alpha C_2}+(\sup_M B)e^{\beta C_2}+C_1
\end{align}
for any $(x,t)\in M\times [0,T)$. Using \eqref{2.6}, \eqref{2.8} and the elliptic $L^p$-estimate yields
\begin{align}
\|u\|_{W^{2,p}(M)}\leq C(\|L u\|_{L^p(M)}+\|u\|_{L^p(M)})\leq C_3
\end{align}
for any $p>1$, where $C_3$ is a constant independent of $t\in [0,T)$. From the Sobolev's embedding theorem, we obtain that $\|u\|_{C^{1,\alpha}(M)}$ is bounded uniformly in $t\in [0,T)$, and hence $\|-S-Ae^{\alpha u}+Be^{-\beta u}\|_{C^\alpha(M)}$ is bounded uniformly in $t\in [0,T)$, too. From \cite[Theorem 2.2.1, pp. 79-80]{[Jost]}, it follows that there exists a constant $C_4>0$ independent of $t\in [0,T)$ such that
\begin{align}
\|u\|_{C^{2,\alpha}(M)}+\left\|\frac{\partial u}{\partial t}\right\|_{C^\alpha(M)}\leq C_4.
\end{align}
Therefore, $u(\cdot,t)\rightarrow u(\cdot, T)$ in $C^{2,\alpha}(M)$, as $t\rightarrow T$, and the short-time existence for parabolic equations implies that the solution $u(x,t)$ can be continued up to some time $T+\epsilon$. Consequently, the solution $u(x,t)$ is global, that is $T=+\infty$. For the uniqueness, let $u_1$ and $u_2$ be two solutions of \eqref{2.1}, then
\begin{align}
\left(\frac{\partial}{\partial t}-L\right)(u_1-u_2)^2=-2&|\nabla(u_1-u_2)|^2-2(u_1-u_2)A(e^{\alpha u_1}-e^{\alpha u_2})\\
&+2(u_1-u_2)B(e^{-\beta u_1}-e^{-\beta u_2})\leq0,\notag
\end{align}
\begin{align}
(u_1-u_2)(x,0)=0.
\end{align}
From Lemma \ref{lemma2.1}, we deduce that $u_1\equiv u_2$.
\qed

\begin{lemma}\label{lemma2.4}
 Let $u(x,t):M\times [0,+\infty)\rightarrow \mathbb{R}$ be the solution of \eqref{2.1}. There exists a constant $C_5>0$ independent of $t$ such that
 \begin{align}
 u(x,t)>-C_5
 \end{align}
 for any $(x,t)\in M\times [0,\infty)$.
\end{lemma}
\proof Set $v=-u+w$, where $w$ is a solution of 
\begin{align}
Lw=S-\bar{S},\ \bar{S}=\text{Vol(M)}^{-1}\int_M S,
\end{align}
where $L=\Delta -\langle d\cdot, \theta\rangle$ with $d^*\theta=0$. Then by \eqref{2.1}, we have
\begin{align}
\frac{\partial v}{\partial t}=Lv+\bar{S}+(Ae^{\alpha w})e^{-\alpha v}-(Be^{-\beta w})e^{\beta v}.
\end{align}
Thus,
\begin{align}\label{2.16}
\frac{\partial v}{\partial t}\leq Lv+\bar{S}+(\sup_M A) e^{\alpha \sup_Mw}e^{-\alpha v},
\end{align}
because $B\geq 0$. Now applying Lemma \ref{lemma2.1} to \eqref{2.16}, we obtain
\begin{align}\label{2.17}
v(x,t)\leq \phi(t),
\end{align}
where $\phi(t)$ is the solution of 
\begin{equation}
\left\{
\begin{aligned}
&\frac{d\phi}{dt}=\bar{S}+(\sup_M A) e^{\alpha \sup_Mw}e^{-\alpha \phi},\\
&\phi(0)=c_0,
\end{aligned}
\right.
\end{equation}
where $c_0$ is a constant with
\begin{align}
c_0>\max\{ \sup_M |-u_0+w|,\ \alpha^{-1}\ln(-\bar{S}^{-1}\sup_M A\cdot e^{\alpha \sup_M w})\}.
\end{align}
By a simple computation, we get
\begin{align}
\phi(t)=\alpha^{-1}\ln{\left(\frac{\sup_M A\cdot e^{\alpha \sup_Mw}+\bar{S}e^{\alpha c_0}}{\bar{S}}e^{\alpha\bar{S}t}-\frac{(\sup_M A) e^{\alpha \sup_Mw}}{\bar{S}}\right)},
\end{align}
which implies that $\phi(t)<C$ for some constant $C>0$ independent of $t\in [0,+\infty)$, because $\bar{S}<0$. Therefore, by \eqref{2.17}, we have
\begin{align}
u(x,t)>\inf_M w-C,
\end{align}
which finishes the proof of this lemma.
\qed

In order to prove the convergence of the flow \eqref{2.1}, we need to establish the uniform $C^0$-estimate for the solution $u(x,t)$.

\begin{lemma}\label{lemma2.3}
 Let $u(x,t):M\times [0,+\infty)\rightarrow \mathbb{R}$ be the solution of \eqref{2.1}. Then
 \begin{align}
 \|u\|_{C^0(M)}\leq C_6 \|u\|_{L^2(M)}+C_6.
 \end{align}
 where $C_6>0$ is a constant independent of $t$.
 \end{lemma}
 \proof Set $v=u+C_5$, where $C_5>0$ is the constant as in Lemma \ref{lemma2.4}, so $v>0$ in $M$. Now multiplying by $v^{a+1}\ (a\geq 0)$ on both sides of \eqref{2.1} and integrating the resulting equation on $M$ yield
\begin{align}
0&=\int_M\left(Lv-\frac{\partial v}{\partial t}-S-Ae^{-\alpha C_5} e^{\alpha v}+Be^{\beta C_5}e^{-\beta v}\right)v^{a+1}\\
&\leq-\frac{4(a+1)}{(a+2)^2}\int_M \left|\nabla v^{\frac{a}{2}+1}\right|^2+\left(C_1+\sup_M|S|+e^{\beta C_5}\sup_M B\right)\int_M v^{a+1},\notag
\end{align}
where we have used the integration by parts, \eqref{2.5}, $d^*\theta=0$, $A\geq 0$ and $v>0$. Hence, we get
\begin{align}
\int_M \left|\nabla v^{\frac{a}{2}+1}\right|^2\leq C\frac{(a+2)^2}{4(a+1)}\int_M v^{a+1}
\end{align}
Using the H\"older's inequality (note $\text{Vol(M)}=1$)
\begin{align}
\int_Mv^{a+1}\leq\frac{a+1}{a+2} \int_Mv^{a+2}+\frac{1}{a+2},
\end{align}
we get
\begin{align}
\int_M \left|\nabla v^{\frac{a}{2}+1}\right|^2\leq C(a+1)\int_M v^{a+2}+C.
\end{align}
In terms of the Sobolev's inequality, we deduce that
\begin{align}\label{1.25}
\left(\int_Mv^{(a+2)\beta}\right)^{\frac{1}{\beta}}&\leq C\left(\int_M\left|\nabla v^{\frac{a}{2}+1}\right|^2+\int_M v^{a+2} \right)\\
&\leq C(a+2)\max\left\{\int_M v^{a+2},1\right\},\notag
\end{align}
where $\beta=\frac{n}{n-2}$. Set $p=a+2$, then it follows from \eqref{1.25} that
\begin{align}
\max\{\|v\|_{L^{p\beta}(M)},1\}\leq C^{\frac{1}{p}}p^{\frac{1}{p}}\max\{\|v\|_{L^p(M)},1\}.
\end{align}
By iteration, we conclude that
\begin{align}
\max\{\|v\|_{L^{p\beta^{k+1}}(M)},1\}&\leq C^{\frac{1}{p\beta^k}}(p\beta^k)^{\frac{1}{p\beta^k}}\max\{\|v\|_{L^{p\beta^k}(M)},1\}\\
&\leq \dots\notag\\
&\leq C^{\sum_{i=0}^k\frac{1}{p\beta^i}}p^{\sum_{i=0}^k\frac{1}{p\beta^i}}\beta^{\sum_{i=0}^k\frac{i}{p\beta^i}}\max\{\|v\|_{L^p(M)},1\}.\notag
\end{align}
Pick $p=2$ and let $k\rightarrow +\infty$, we obtain 
\begin{align}
\|v\|_{C^0(M)}\leq \max\{\|v\|_{C^0(M)},1\}\leq C\max\{\|v\|_{L^2(M)},1\}.
\end{align}
Since $v=u+C_5$, then by the triangle inequality we have
\begin{align}
\|u\|_{C^0(M)}\leq C\|u\|_{L^2(M)}+C.
\end{align}
 \qed

\begin{lemma}\label{lemma2.6}
Let $u(x,t):M\times [0,+\infty)\rightarrow \mathbb{R}$ be the solution of \eqref{2.1}. Then
\begin{align}
\|u\|_{L^2(M)}<C_7,
\end{align} 
 where $C_7>0$ is a constant independent of $t$.
\end{lemma}

\proof We prove it by contradiction. Suppose that $\|u\|_{L^2(M)}$ is not uniformly bounded with respect to $t\in[0,+\infty)$, then there exists a sequence $t_i\rightarrow +\infty$ such that 
\begin{align}
\lim_{i\rightarrow +\infty}\|u(\cdot,t_i)\|_{L^2(M)}=+\infty,\ \frac{d}{dt}\vert_{t=t_i}\|u(\cdot,t)\|_{L^2(M)}\geq 0.\label{3.26.}
\end{align}
Set 
\begin{align}\label{2.23}
u_i=u(\cdot,t_i),\ l_i=\|u(\cdot,t_i)\|_{L^2(M)},\ w_i=\frac{u_i}{l_i}.
\end{align}
From \eqref{2.23} and Lemma \ref{lemma2.3}, we deduce that
\begin{align}\label{2.24}
\|w_i\|_{L^2(M)}=1,\ \|w_i\|_{C^0(M)}\leq C,
\end{align}
where $C>0$ is a constant independent of $i$. By \eqref{2.1} and H\"older's inequality, we perform the following computation
\begin{align}
\int_M|\nabla u_i|^2&=-\int_Mu_i Lu_i\label{3.28}\\
&=\int_M u_i\left(-S-A e^{\alpha u_i}+B e^{-\beta u_i}-\frac{\partial u}{\partial t}\vert_{t=t_i}\right)\notag\\
&\leq \left(\sup_M |S|+\sup_{M\times[0,+\infty)}\left|\frac{\partial u}{\partial t}\right|\right)\int_M |u_i|+\int_M\left(\frac{A}{\alpha}+\frac{B}{\beta}\right)\notag\\
&\leq C\|u_i\|_{L^2(M)}+\int_M\left(\frac{A}{\alpha}+\frac{B}{\beta}\right),\notag
\end{align}
where $C>0$ is a constant independent of $i$, and we have used that $d^*\theta=0$ and $xe^x\geq-1\ (\forall x\in \mathbb{R})$. Dividing \eqref{3.28} by $l_i^2$ gives
\begin{align}
\int_M|\nabla w_i|^2\leq \frac{C_{9}}{l_i}+\frac{\int_M\left(\frac{A}{\alpha}+\frac{B}{\beta}\right)}{l_i^2},
\end{align}
hence,
\begin{align}\label{2.27}
\lim_{\ i\rightarrow +\infty}\int_M|\nabla w_i|^2=0.
\end{align}
By \eqref{2.24} and \eqref{2.27}, we get that $\{w_i\}_{i=1}^\infty$ is bounded in $W^{1,2}(M)$. According to the Sobolve's embedding theorem $W^{1,2}(M)\subset\subset L^2(M)$, there is a subsequence of $\{w_i\}_{i=1}^\infty$, also denoted by $\{w_i\}_{i=1}^\infty$, and $w_{\infty}\in W^{1,2}(M)$ such that
\begin{align}
&w_i\rightarrow w_{\infty}\ \text{in\ } L^2(M),\\
&w_i\rightharpoonup w_{\infty}\ \text{weakly\ in\ } W^{1,2}(M).\label{3.32}
\end{align}
From\eqref{2.24} and \eqref{2.27}-\eqref{3.32}, we deduce that
\begin{align}
w_\infty=C^*\not=0\ \text{almost\ everywhere\ in\ } M,
\end{align} 
where $C^*$ is a constant. From \eqref{2.1}, \eqref{2.5} and Lemma \ref{lemma2.4}, it follows that
\begin{align}
\int_MAe^{\alpha u_i}&=\int_M\left(Lu_i-S+Be^{-\beta u_i}-\frac{\partial u}{\partial t}|_{t=t_i}\right)\\
&\leq -\int_MS+\int_M Be^{\beta C_5}+C_1\text{Vol(M)}.\notag
\end{align}
Hence,
\begin{align}\label{3.35}
\alpha\int_MAw_id\mu_\eta&\leq\frac{\int_MAe^{\alpha l_iw_i}}{l_i}\\
&\leq\frac{-\int_MS+\int_M Be^{\beta C_5}+C_1\text{Vol(M)}}{l_i}\rightarrow 0, \text{as}\ i\rightarrow +\infty.\notag
\end{align}
In terms of \eqref{3.32}, we have
\begin{align}
\lim_{i\rightarrow +\infty}\int_MAw_i=C^*\int_MA.\label{3.36}
\end{align}
Combining \eqref{3.35} and \eqref{3.36} yields $C^*<0$, since $A\geq 0\ (\not\equiv 0)$ and $C^*\not=0$. Using \eqref{2.1}, \eqref{3.32} and the fact $xe^x\geq -1\ (\forall \ x\in \mathbb{R})$, we obtain 
\begin{align}
\frac{d}{dt}|_{t=t_i}\|u\|_{L^2(M)}&=\frac{1}{l_i}\int_Mu_i\frac{\partial u}{\partial t}|_{t=t_i}\\
&=\frac{1}{l_i}\int_Mu_i\left(L u_i-S-Ae^{\alpha u_i} +Be^{-\beta u_i}\right)\notag\\
&\leq-\frac{1}{l_i}\int_M|\nabla u_i|^2-\int_Mw_iS+\frac{\int_MA}{\alpha l_i}+\frac{\int_M B}{\beta l_i}\notag\\
&\leq-\int_Mw_iS+\frac{\int_MA}{\alpha l_i}+\frac{\int_M B}{\beta l_i}\rightarrow -C^*\int_MS,\notag
\end{align}
as $i\rightarrow +\infty$, which is a contradiction with \eqref{3.26.}, because $\int_MS<0$ and $C^*<0$.
\qed

 \begin{remark}
The Lemma \ref{lemma2.3} can also be proved by a similar technique in the proof of \cite[Lemma 3.3]{[Yu]} without using Lemma \ref{lemma2.4}. However, in our case, Lemma \ref{lemma2.4} is necessary for the proof of Lemma \ref{lemma2.6}.
 \end{remark}

Combining Lemma \ref{lemma2.3} and Lemma \ref{lemma2.6} gives
\begin{lemma}\label{lemma2.7}
 Let $u(x,t):M\times [0,+\infty)\rightarrow \mathbb{R}$ be the solution of \eqref{2.1}. Then
 \begin{align}
 \|u\|_{C^0(M)}\leq C,
 \end{align}
 where $C>0$ is a constant independent of $t$.
 \end{lemma}
 
Now using this uniform $C^0$-estimate for the solution of \eqref{2.1}, we prove the following two convergence results for this flow.

 \begin{theorem}\label{theorem2.8}
 Let $u(x,t):M\times [0,+\infty)\rightarrow \mathbb{R}$ be the solution of \eqref{2.1}. If either $A>0$ or $B>0$, then there exists a sequence $t_i\rightarrow +\infty$ such that $u(\cdot,t_i)$ converges to $u_\infty$ in $C^{1,\alpha}(M)$, where $\alpha\in (0,1)$, and $u_\infty$ is the unique smooth solution of 
 \begin{align}\label{2.37}
 \Delta u_\infty-\langle du_\infty, \theta\rangle-S-A e^{\alpha u_\infty}+B e^{-\beta u_\infty}=0.
 \end{align}
 \end{theorem}
\proof By \eqref{2.5} and Lemma \ref{lemma2.7}, it follows from \eqref{2.1} that
\begin{align}
\|Lu\|_{C^0(M)}\leq C,
\end{align}
where $C$ is a constant independent of $t\in[0,+\infty)$. Hence, by elliptic $L^p$-estimate, we obtain that $\|u\|_{W^{2,p}(M)}$ is uniformly bounded with respect to $t\in[0,+\infty)$ for any $p>1$. In terms of the Sobolev's embedding theorem, there exists a sequence $t_i\rightarrow +\infty$ and $u_\infty\in C^{1,\alpha}(M)$ such that
\begin{align}\label{2.39}
u(\cdot,t_i)\rightarrow u_\infty \text{\ in\ } C^{1,\alpha}(M).
\end{align}
According to \eqref{2.3} and Lemma \ref{lemma2.7}, we obtain that
\begin{align}\label{2.38}
\left(\frac{\partial}{\partial t}-L\right)\left|\frac{\partial u}{\partial t}\right|^2\leq -C\left|\frac{\partial u}{\partial t}\right|^2,
\end{align}
where $C=2\left(\alpha  e^{\alpha \|u\|_{C^0(M)}} \inf_M A+\beta e^{\beta \|u\|_{C^0(M)}} \inf_M B \right)>0$ is a constant independent of $t$, since either $A>0$ or $B>0$ on $M$. Now applying Lemma \ref{lemma2.1} to \eqref{2.38} yields
\begin{align}
\left|\frac{\partial u}{\partial t}\right|^2\leq e^{-Ct}\sup_M\left|L u_0-S-A e^{\alpha u_0}+B e^{-\beta u_0}\right|^2,
\end{align}
thus,
\begin{align}\label{2.42}
\lim_{t\rightarrow +\infty}\max_M\left|\frac{\partial u}{\partial t}\right|^2=0.
\end{align}
Combining \eqref{2.39} with \eqref{2.42}, we deduce that $u_\infty$ is a weak solution of \eqref{2.37}. By the regularity results for elliptic equations, we can conclude that $u_\infty\in C^\infty(M)$. For the uniqueness, let $u_1, u_2$ be two smooth solutions of \eqref{2.37}, and set $v=u_1-u_2$, then
\begin{align}
\Delta v-\langle dv, \theta\rangle+A e^{\alpha u_2}\left(1-e^{\alpha v}\right)+B e^{-\beta u_1}\left(1-e^{\beta v}\right)=0.
\end{align}
Using the maximum principle, it is clear that $v\equiv 0$, i.e., $u_1=u_2$.
\qed
\begin{theorem}\label{theorem2.9}
 Let $u(x,t):M\times [0,+\infty)\rightarrow \mathbb{R}$ be the solution of \eqref{2.1} with $\theta\equiv 0$. Then there exists a sequence $t_i\rightarrow +\infty$ such that $u(\cdot,t_i)$ converges to $u_\infty$ in $C^{1,\alpha}(M)$, where $\alpha\in (0,1)$, and $u_\infty$ is the unique smooth solution of 
 \begin{align}\label{2.37}
 \Delta u_\infty-S-A e^{\alpha u_\infty}+B e^{-\beta u_\infty}=0.
 \end{align}
 \end{theorem}
 \proof
 Since
 \begin{align}
 \int_M\left|\frac{\partial u}{\partial t}\right|^2&=\int_M \left(\Delta u-S-A e^{\alpha u}+B e^{-\beta u}\right)^2\\
 &=-\frac{d}{dt}\int_M\left(\frac{1}{2}|\nabla u|^2+Su+\frac{A}{\alpha}e^{\alpha u}+\frac{B}{\beta}e^{-\beta u}\right),\notag
 \end{align}
 then from Lemma \ref{lemma2.7} we have
 \begin{align}
\int_{0}^{t} \int_M \left|\frac{\partial u}{\partial t}\right|^2&\leq-\int_M\left(Su+\frac{A}{\alpha}e^{\alpha u}+\frac{B}{\beta}e^{-\beta u}\right)\\
&\ \ \ \ \ \ +\int_M\left(\frac{1}{2}|\nabla u_0|^2+Su_0+\frac{A}{\alpha}e^{\alpha u_0}+\frac{B}{\beta}e^{-\beta u_0}\right)\notag\\
&\leq C\notag,
\end{align}
where $C>0$ is a constant independent of $t$. Hence, 
\begin{align}
\int_{0}^{\infty} \int_M \left|\frac{\partial u}{\partial t}\right|^2<\infty.
\end{align}
So there is a sequence $t_i\rightarrow +\infty$ such that 
\begin{align}
\int_M \left|\frac{\partial u}{\partial t}\right|^2(\cdot, t_i)\rightarrow 0.
\end{align}
By a similar argument in Theorem \ref{theorem2.8}, we can also obtain a subsequence $t_i\rightarrow +\infty$ such that $u(\cdot, t_i)\rightarrow u_\infty$ in $C^{1,\alpha}(M)$, where $u_\infty$ is the unique $C^\infty(M)$ solution of \eqref{2.37}.
\qed

For the case that $A\geq 0\ (\not\equiv 0),\ B\geq 0$ and $d^*\theta=0$, the above flow approach seems difficult to handle it. Here we apply the method of upper and lower solutions to deal with this case.

\proof[\textbf{Proof of Theorem \ref{theorem1.1}}] Let $v_1, v_2\in C^\infty(M)$ be solutions of 
\begin{align}
Lv_1=S-\overline{S},\ Lv_2=A-\overline{A},
\end{align}
respectively, where $\bar{S}=\text{Vol(M)}^{-1}\int_M S<0$ and $\bar{A}=\text{Vol(M)}^{-1}\int_M A>0$, and $L=\Delta -\langle d\cdot, \theta\rangle$ with $d^*\theta=0$. Set $u_+=v_1+av_2+b$, where $a,b>0$ are two constants that we will choose later, and let $c_1=\inf_Mv_1$ and $c_2=\inf_M v_2$, then
\begin{align}
&\ \ \ \ \ \ Lu_+-S-A e^{\alpha u_+}+B e^{-\beta u_+}\\
&\leq-\bar{S}-a\bar{A}+A\left(a-e^{\alpha(c_1+ac_2+b)}\right)+e^{-\beta(c_1+ac_2+b)}\sup_MB<0,\notag
\end{align}
provided that $a, b>0$ satisfy
\begin{equation}
\left\{
\begin{aligned}
&a>-\frac{\bar{S}}{\bar{A}},\\
&b>\max\left\{\alpha^{-1}\ln a-c_1-ac_2,\ -c_1-ac_2-\beta^{-1}\ln{\frac{\bar{S}+a\bar{A}}{\sup_M B}}\right\}.\\
\end{aligned}
\right.
\end{equation}
Define $u_-=v_1-m$, where $m$ is a positive number, then
\begin{align}
u_-<u_+,
\end{align}
\begin{align}
Lu_--S-A e^{\alpha u_-}+B e^{-\beta u_-}\geq -\bar{S}-e^{\alpha (\sup_M v_1-m) }\sup_M A>0,
\end{align}
provided that $m$ satisfies
\begin{align}
m>\max\left\{\sup_M v_1-\alpha^{-1}\ln \frac{-\bar{S}}{\sup_M A},\ -ac_2-b\right\}.
\end{align}
Using the method of upper and lower solutions (e.g., \cite[Lemma 3.4]{[BW]}), it is easy to see that there exists a smooth solution $u$ of \eqref{1.3} with $u_-\leq u\leq u_+$. For the uniqueness, it follows from the maximum princple.

\qed
\section{Proof of Proposition \ref{prop1.2}}
In this section, we will prove Proposition \ref{prop1.2}. To this end, we need the following lemmas.
\begin{lemma}\label{lemma3.2}
There exists a constant $C>0$ such that for any $u\in C^\infty(M)$ satisfying \eqref{1.3} we have
\begin{align}
u(x)>-C,\  \forall \ x\in M,
\end{align}
\end{lemma}
\proof Set $\phi_t=w-t$, where $t\in \mathbb{R}$ and $w$ is a solution of 
\begin{align}
Lw=S-\bar{S},\ \bar{S}=\text{Vol(M)}^{-1}\int_M S.
\end{align}
Then for any $t>t_0=-\alpha^{-1}\ln\left\{ -(\sup_M A)^{-1}e^{-\alpha\|w\|_{C^0(M)}}\bar{S}\right\}$, we have
\begin{align}\label{3.5}
 L \phi_t-S-A e^{\alpha \phi_t}+B e^{-\beta \phi_t}\geq-\bar{S}-e^{\alpha \|w\|_{C^0(M)}}e^{-\alpha t}\sup_MA>0.
\end{align}
Now we claim that $u\geq \phi_{t_0}$ on $M$. Otherwise, there exists a constant $t_1\in (t_0,+\infty)$ such that $u\geq \phi_{t_1}$ on $M$ and $u(x_0)=\phi_{t_1}(x_0)$. Hence, by the maximum principle, at $x_0$, we have
\begin{align}
&\Delta \phi_{t_1}-\langle d\phi_{t_1}, \theta\rangle-S-A e^{\alpha \phi_{t_1}}+B e^{-\beta \phi_{t_1}}\\
&=\Delta (\phi_{t_1}-f)-\langle d(\phi_{t_1}-f), \theta\rangle+\left(\Delta f-\langle df, \theta\rangle-S-A e^{\alpha f}+B e^{-\beta f}\right)\notag\\
&\ \ \ \ \ \ +A(e^{\alpha f}-e^{\alpha \phi_{t_1}})+B( e^{-\beta \phi_{t_1}}- e^{-\beta f})\leq 0,\notag
\end{align}
which leads to a contradiction with \eqref{3.5}. Therefore, $u\geq \min_M\phi_{t_0}$. 
\qed

\begin{lemma}\label{lemma3.2.}
There exists a constant $C>0$ such that for any $u\in C^\infty(M)$ satisfying \eqref{1.3} we have
\begin{align}
\|u\|_{L^2(M)}<C.
\end{align}
\end{lemma}
\proof Since
\begin{align}
\int_M |u|^2=\int_{\{-C<u\leq 0\}}|u|^2+\int_{\{u>0\}}(u^+)^2\leq C^2\text{Vol(M)}+\int_{M}(u^+)^2,
\end{align}
where $u^+=\max\{u, 0\}$, so it is sufficient to prove $\|u^+\|_{L^2(M)}\leq C$. Pick a test function $\phi=u^+\in W^{1,2}(M)$ in \eqref{1.3}, we obtain
\begin{align}\label{3.9}
\int_M\left(-|\nabla u^+|^2-\langle du, \theta\rangle u^+- Su^+- A e^{\alpha u}u^++B e^{-\beta u}u^+\right)=0.
\end{align}
Using $d^*\theta=0$, we get
\begin{align}\label{3.10}
\int_M\langle du, \theta\rangle u^+=\frac{1}{2}\int_M\langle d(u^+)^2, \theta\rangle=\frac{1}{2}\int_M\langle (u^+)^2, d^*\theta\rangle=0.
\end{align}
In terms of Young's inequality and Lemma \ref{lemma3.2}, we deduce that
\begin{align}
-\int_M Su^+\leq (4\delta_1)^{-1}\int_M S^2+\delta_1\int_M(u^+)^2,\label{3.11}\\
\left|\int_MB e^{-\beta u}u^+\right|\leq e^{\beta C}\sup_M B\left(\delta_2\int_M (u^+)^2+(4\delta_2)^{-1}\text{Vol(M)} \right).\label{3.12}
\end{align}
As $A>0$ and $e^{2t}\geq t^3$, $t^4\geq t^2-1$ for $t\in \mathbb{R}$, we have
\begin{align}\label{3.13}
\int_MA e^{\alpha u}u^+\geq \frac{\alpha^3}{8}\eta\int_M(u^+)^4\geq\frac{\alpha^3}{8}\eta\left(\int_M(u^+)^2-\text{Vol(M)}\right),
\end{align}
where $\eta=\inf_M A>0$.
Substituting \eqref{3.10}-\eqref{3.13} into \eqref{3.9} yields
\begin{align}
\left(\frac{\alpha^3\eta}{8}-\delta_1- \tilde{C}\delta_2\right)\int_M(u^+)^2\leq \frac{1}{4\delta_1}\int_M S^2+\left(\frac{\tilde{C}}{4\delta_2}+\frac{\alpha^3\eta}{8}\right)\text{Vol(M)},
\end{align}
where $\tilde{C}=e^{\beta C}\cdot \sup_M B$. Taking $\delta_1=\frac{\alpha^3\eta}{24}$ and $\delta_2=\frac{\alpha^3\eta}{24\tilde{C}}$ gives
\begin{align}
\int_M(u^+)^2\leq \frac{144}{\alpha^6}\eta^{-2}\int_M S^2+\left(\tilde{C}\frac{144}{\alpha^6}\eta^{-2}+3\right)\text{Vol(M)}.
\end{align}
\qed

\begin{remark}
The analytic techniques implemented in the proof of Lemma \ref{lemma3.2.} cannot be applied to prove Lemma \ref{lemma2.6}, because in the proof of Lemma \ref{lemma3.2.} we need to use the assumption $A>0$ which implies $\inf_M A>0$.
\end{remark}

\begin{lemma}\label{lemma3.4}
There exists a constant $C>0$ such that for any $u\in C^\infty(M)$ satisfying \eqref{1.3} we have
\begin{align}
\|u\|_{C^0(M)}<C\|u\|_{L^2(M)}+C.
\end{align}
\end{lemma}
\proof The proof of this lemma is similar to that of Lemma \ref{lemma2.3}, so we omit it here.

\proof[\textbf{Proof of Proposition \ref{prop1.2}}] From Lemma \ref{lemma3.2.} and Lemma \ref{lemma3.4}, we have $\|u\|_{L^p(M)}\leq C$ and $ \|S+A e^{\alpha u}-B e^{-\beta u}\|_{L^p(M)}\leq C$ for any $p>1$. According to the elliptic $L^p$-estimate, we get $\|u\|_{W^{2,p}(M)}\leq C$ for any $p>1$. By the Sobolev's embedding theorem, it follows that $\|u\|_{C^{1,\alpha}(M)}\leq C$ for any $\alpha\in (0, 1)$. Now using the Schauder estimate for elliptic equations, it is easy to see that \eqref{3.2} holds.
\qed

\bigskip
\bigskip

Weike Yu

School of Mathematics and Statistics

Nanjing University of Science and Technology

Nanjing, 210094, Jiangsu, P. R. China

wkyu2018@outlook.com

\bigskip

\end{document}